\documentclass[reqno,11pt]{amsart}

\usepackage{color}

\newtheorem{theorem}{Theorem}

\theoremstyle{definition}
\newtheorem{assumption}[theorem]{Assumption}

\theoremstyle{remark}
\newtheorem{remark}[theorem]{Remark}

\makeatletter 
 \def\dashint{%
 \operatorname%
 {\,\,\text{\bf--}\kern-.98em\DOTSI\intop\ilimits@\!\!}}
\def\dashnorm{\,\,\text{\bf--}\kern-.5em\|}
\def\ninf{\qopname\relax\@empty{inf\phantom{p}\!\!\!}}

 \makeatother

\newcommand\bB{\mathbb{B}}
\newcommand\bC{\mathbb{C}}

\newcommand\bR{\mathbb{R}}

\newcommand\cF{\mathcal{F}}

\newcommand{\loc}{{\rm loc}\,}

\begin{document}

\title[Remark on R\"ockner and Zhao]
{Remarks on the post of R\"ockner and Zhao
about strong solutions of stochastic equations
with critical singularities
at arXiv:2103.05803v4}
\author{N.V. Krylov}

\email{nkrylov@umn.edu}
\address{School of Mathematics, University of Minnesota, Minneapolis, MN, 55455}
 
 

\begin{abstract} 
 We discuss various relations between
part of the results in arXiv:2103.05803v4
and recent author's results.
\end{abstract}

 \maketitle

                                                  \label{section 3.11.1}

Let $\bR^{d}$ be a $d-$dimensional Euclidean space of points
$x=(x^{1},...,x^{d})$ with $d\geq3$. Let $(\Omega,\cF,P)$ be a 
complete probability space,
carrying a $d $-dimensional Wiener process
$w_{t}$ and let $b(t,x)$ be a Borel
measurable $\bR^{d}$-valued function on $\bR^{d+1}
=\{(t,x):t\in\bR, x\in \bR^{d}\}$. 
 In the really mind opening paper
of   M.~R\"ockner and G.~Zhao \cite{RZ_21},
  the authors prove among many other things
 the following, which is part of their
 Theorem 1.1. 
Take $T\in(0,\infty)$ and set
$\bR^{d}_{T}=[0,T]\times\bR^{d}$.

\begin{theorem}
                           \label{theorem 7.19.1}

Assume that $b$ satisfies one of the following two conditions:

(a) $b\in C\big([0,T];L_{d}(\bR^{d})\big)$;

(b) $b\in L_{p_{1},q_{1}}(\bR^{d}_{T}):=L_{q_{1}}\big([0,T];L_{p_{1}}(\bR^{d})\big)$ with $p_{1}\in(d,\infty),q_{1}\in(2,\infty)$ and $d/p_{1}+2/q_{1}=1$.

Fix $x\in\bR^{d}$. Then the equation
\begin{equation}
                               \label{7.19.1}
x_{t}=x+\int_{0}^{t}b(s,x_{s})\,ds+w_{t}
\end{equation}
admits a  strong solution 
such that,
for {\em any\/} $p,q $ satisfying
\begin{equation}
                          \label{7.19.4}
p,q\in(1,\infty), \quad
\frac{d}{p}+\frac{2}{q}<2
\end{equation}
and any Borel nonnegative $f$ on $\bR^{d+1}$,
we have
\begin{equation}
                               \label{7.19.2}
E\int_{0}^{T}f(t,x_{t})\,dt\leq
N\sup_{y\in\bR^{d}}\|fI_{B_{1}(y)}\|_{
L_{p ,q }(\bR^{d}_{T})},
\end{equation}
where  $N$
is a constant independent of $f$ and $x$.
Furthermore, there is only one 
strong solution
of \eqref{7.19.1} such that \eqref{7.19.2} holds
for {\em any\/} $p,q$ satisfying \eqref{7.19.4}
 (with $q_{1}=\infty$ in case(a)) and any Borel nonnegative $f$
with a constant  $N$
 independent of $f$.

\end{theorem}

The purpose of this note is to compare
Theorem \ref{theorem 7.19.1} with the following
simplified version of Theorem 2.4 from \cite{Arx2}
under the assumption that follows, in which $\bB_{\rho}$
is the collection of balls of radius $\rho$
in $\bR^{d}$, and $\bC_{\rho}$ is the collection
of cylinders $C=(s,s+\rho^{2})\times B$, where
$B\in\bB_{\rho}$,
$$
\dashnorm f\|_{L_{p}(B)}^{p}=\frac{1}{\text{Vol}(B)}
\int_{B}|f|^{p}\,dx,\quad\dashnorm f\|^{q}_{L_{p,q}(C)}
=\rho^{-2}\int_{s}^{s+\rho^{2}}\dashnorm f(t,\cdot)\|_{L_{p}(B)}^{q}\,dt
$$
$\rho_{b}\in(0,1)$, 
$p_{b}\in(d/2,d]$, $p_{b}>2$,  
\begin{equation}
                         \label{5.31.2}
 p_{0}\in (2\vee(d/2), p_{b}),\quad q_{0}\in(1,\infty),
\quad \frac{d}{p_{0}}+\frac{2}{q_{0}}<2.
\end{equation}
\begin{assumption}
                          \label{assumption 8.1.1}
We have
$b = b_{M} +b_{B}  $, where $b_{M}$ and $b_{B}$ are   Borel functions  such that there exists 
a constant $\hat b_{M}<\infty$ for which
\begin{equation}
                             \label{6.12.3}
\dashnorm b_{M}(t,\cdot)\|_{L_{p_{b} }(B)}\leq \hat b_{M}  \rho^{-1},
\end{equation}
as long as $t\in[0,T]$, $B\in\bB_{\rho}$, and $\rho\leq \rho_{b}$,
and  
$$
\beta_{b}(t)=\sup_{s}\int_{s}^{s+t} \bar b^{2}_{B} (r)   \,dr<\infty,\quad \bar b_{B}(r):=
\sup_{x\in\bR^{d}}|b_{B} (r, x)|.
$$
 
\end{assumption}

Here are our main results.

\begin{theorem}
                        \label{theorem 6.12.2}
There exists a  constant  $ \alpha_{b}>0$,
depending only on $d, \delta$,   $
p_{b} $,   $p_{0},q_{0}$,    such that, if 
\begin{equation}
                           \label{6.12.4}
 \hat b_{M} \leq \alpha_{b},
\end{equation}
then equation \eqref{7.19.1} has a strong
solution  possessing the property 

(a)  for
 any $T\in(0,\infty)$, $m=1,2,...$,
there exists a constant $N$
such that for any Borel nonnegative $f$
\begin{equation}
                                \label{5.31.1}
E \Big(\int_{0}^{T}f( s,x_{s})\,ds\Big)^{m}\leq 
N \|f\|^{m}_{L_{p_{0},q_{0} }}.
\end{equation}
 
Furthermore (conditional uniqueness), if there exists a  solution $y_{s}$ of \eqref{7.19.1} with the same initial condition having the property 

($\text{a}\,'$) for each $T\in(0,\infty)$
there is a constant $N$ such that for any Borel nonnegative $f$
\begin{equation}
                                \label{7.19.01}
E  \int_{0}^{T}f( s,y_{s})\,ds \leq 
N \|f\| _{L_{p_{0},q_{0}}},
\end{equation}
then $x_{\cdot}=y_{\cdot}$ (a.s.).
\end{theorem}

The statements about uniqueness in Theorems 
\ref{theorem 7.19.1} and \ref{theorem 6.12.2} are somewhat conditional,
however, in \cite{Arx2} one can find a substantial
amount of cases when any solution is strong and unique.

\begin{remark}
                   \label{remark 1.5.1}

 Let us compare conditions
\eqref{6.12.4} and (b) in Theorem 
\ref{theorem 7.19.1}. Assume that
(b) in Theorem 
\ref{theorem 7.19.1} is satisfied. For simplicity we fix $T\in(0,\infty)$
and suppose that $b(t,x) =0$ for $t\not \in [0,T]$.

{\em Case $p>d$ (and $q>2$)\/}. Set $p_{b}=p_{D\sigma}=d$ ($\sigma$ is the unit matrix).
Then take a  constant $\hat N>0$
and let
$$
\lambda(t)=\hat N\Big(\int_{\bR^{d}}
|b (t,x)|^{p  }\,dx\Big)^{1/(p  -d )}.
$$
Also define  
$$
b_{M} (t, x)=b  (t, x)I_{|b  (t,x)|\geq \lambda(t)}
$$
and observe that  for $B\in\bB_{\rho}$  we have
$$
\dashint_{B }|b_{M} (t, x)|^{d}\,dx
\leq \lambda^{d-p  }(t)
\dashint_{B }|b  (t, x)|^{p  }\,dx
\leq N(d)\hat N^{d-p  }\rho^{-d}.
$$
 
Here $N(d)\hat N^{d-p  }$ is as small
as we like if   $\hat N$
large enough. Furthermore, for 
$b_{B} =b -b_{M} $ it holds that
 $|b_{B} |\leq \lambda $ and
$$
\int_{0}^{T}\lambda^{2}(t)\,dt
=\hat N^{2}\int_{0}^{T}\Big(\int_{\bR^{d}}|b (t,x)|^{p  }
\,dx\Big)^{q  /p  }\,dt<\infty.
$$
This shows that $b$ satisfies Assumption
\ref{assumption 8.1.1} with $\hat b_{M}$ as small as we like.  Thus,   Theorem \ref{theorem 6.12.2}
is applicable and yields a strong solution
conditionally unique.

It is shown in \cite{Arx2} that the solution is, actually, unique
(unconditionally)   if $p\geq d+1$.

If $p\in(d,d+1)$, for a {\em strong\/} solution
$y_{\cdot}$ 
to coincide with the strong one $x_{\cdot}$, property (b), with $y_{\cdot}$ in place of $x_{\cdot}$,
is required in \cite{RZ_21}, where
$ p, q$ are {\em any\/} numbers satisfying \eqref{7.19.4}.
According to Theorem \ref{theorem 6.12.2}
we only need ($\text{a}'$) with {\em some\/} $p,q$ satisfying \eqref{5.31.2} and do not need $y_{\cdot}$ to be a {\em strong\/} solution (it might be that this
discrepancy is caused by a sloppy use of the quantifiers). On the other hand, although the condition on the range of $p,q$
from \cite{RZ_21} is stronger (in uniqueness),   the right-hand
side of \eqref{7.19.01} is smaller than 
the right-hand side of \eqref{7.19.1} 
for $p=p_{0},q=q_{0}$.

{\em Case $p=d$ and $b\in C([0,T],L_{d})$\/}. In that case
\begin{equation}
                                  \label{6.3.1}
 \lim_{r\downarrow 0}\sup_{t\in[0,T]}
\sup_{B\in \bB_{r}}\|b(t,\cdot) \|_{L_{d}(B)}=0.
\end{equation} 

Then, one can take $p_{b} =d$, choose any
$  p_{0} ,q_{0}  $ satisfying \eqref{5.31.2} and set $b_{M}=b_{|b|>1}$,
which guarantees the arbitrary smallness of $\hat b_{M}$ for small enough $\rho_{b}$. Also, obviously,
$$
\int_{0}^{T}\|b(t,\cdot)I_{|b(t,\cdot)|\leq 1}\|^{2}_{L_{\infty}}\,dt<\infty,
$$
  Therefore, this case is also covered by 
Theorem \ref{theorem 6.12.2} and the implications
of that are the same as in the  case of $p\in(d,d+1)$
(a result comparable
to the one in \cite{RZ_21}).

 We see that, actually, condition that
$b\in C([0,T],L_{d})$ can be replaced with $\leq
\varepsilon$ in
\eqref{6.3.1} in place of $=0$, for $\varepsilon>0$ small enough, which holds, for instance,
if the norms $\|b(t,\cdot)\|_{L_{d}}$ are uniformly
sufficiently small, that is imposed as one of alternative conditions
in \cite{RZ_21}.

Another, border case in Theorem \ref{theorem 7.19.1} (b)
when $p=\infty,q=2$, not covered in \cite{RZ_21},
also can be treated by our methods. We take
$b_{M}=0$, $p_{b}=d$,  
any $p_{0}\in (2\vee (d/2),d),q_{0}\in(1,\infty)$,
satisfying \eqref{5.31.2},  and observe that
$\beta_{b}(\infty)<\infty$,
Theorem \ref{theorem 6.12.2}
 is applicable
and yields a strong solution  which is unique
in the set of all solutions possessing property 
($\text{a}'$).

\end{remark}
 
\begin{remark}
                          \label{remark 8.1.3}
It turns out that, generally,
condition \eqref{6.12.4} with $d>p_{b}\geq d-1$ does not imply that
$b\in L_{p_{b}+ \varepsilon,\loc}$, no matter how small $\varepsilon>0$ is. Therefore, such situations
are way beyond the scope of Theorem \ref{theorem 7.19.1} and, since the case (b) is called critical,
we can call our case "supercritical".

Here is
an example. Take $r_{n}>0$, $n=1,2,...$, such that
the  sum of $\rho_{n}:=r_{n}^{d-p_{b}}$ is $1/2$, let $e_{1}$ be the first
basis vector, and set $b(x)=|x|^{-1}
I_{|x|<1}$, $x_{0}=1$,
$$
x_{n}=1-  2\sum_{1}^{n}r_{i}^{d-p_{b}},\quad 
c_{n}=(1/2)(x_{n}+x_{n-1})
$$
$$
  b_{n}(x)=r_{n}^{-1}b\big(r_{n}^{-1}
(x-c_{n}e_{1})\big),\quad b=\sum _{1}^{\infty}b_{n}.
$$
Since $r_{n}\leq 1$ and $d-p_{b}\leq 1$,  the supports of $b_{n}$'s are disjoint and
for $p>0$
$$
\int_{|x|\leq 1}b^{p}\,dx=\sum _{1}^{\infty}\int_{\bR^{d}}b_{n}^{p}\,dx=N(d,p)\sum_{1}^{\infty}r_{n}^{d-p}.
$$
According to this we take the $r_{n}$'s so that
the last sum diverges for any $p>p_{b}$.
Then observe that for any $n\geq 1$ and any ball $B$
of radius $\rho$
$$
 \int_{B } |b_{n}|^{p_{b}}dx \leq N(d) \rho^{d-p_{b}} .
$$
Also, if the intersection of $B$ with $\bigcup B_{r_{n}}(c_{n})$
is nonempty, the intersection
 consists of some $B_{r_{i}}(c_{i})$, $i=i_{0},...,i_{1}$, and $B\cap B_{r_{i_{0}-1}}(c_{i_{0}-1})$ if $i_{0}\ne 0  $ and 
$B\cap B_{r_{i_{1}+1}}(c_{i_{1}+1})$. 
In this situation
$$
 \sum_{i=i_{0} }^{i_{1} }
\rho_{i} \leq 2\rho,
$$
and therefore,
$$
 \int_{B } |b |^{p_{b}}\,dx=N(d)\sum_{i=i_{0}}^{i_{1}}r_{i}^{d-p_{b}}+\int_{B } |b_{i_{0}-1} |^{p_{b}}dx
+\int_{B } |b_{i_{1}+1} |^{p_{b}}dx \leq N(d)(\rho
+\rho^{d-p_{b}}),
$$
where the last term is less than $N(d)\rho^{d-p_{b}}$
for $\rho\leq 1$ and this yields  just a different form of \eqref{6.12.4}.

\end{remark}

 \begin{remark}
                          \label{remark 8.1.7}
As a conclusion we emphasize that
the assumptions of Theorem \ref{theorem 6.12.2}
are much weaker than that of Theorem \ref{theorem 7.19.1}. However, the goals of \cite{RZ_21}
go way beyond of just proving Theorem \ref{theorem 7.19.1}, with one of the main goals being proving
weak differentiability of strong solutions with respect to the initial starting point.   On the other hand, in contrast to
\cite{RZ_21},  the diffusion matrix in \cite{Arx2} is not constant
or even continuous.

\end{remark}

\end{document}